\documentclass[fleqn,10pt]{wlscirep}
\usepackage{color}
\title{Noise Response Data Reveal \ahmet{Novel Controllability Gramian for} Nonlinear Network Dynamics \template{(20 words or less)}}

\author[1,*]{Kenji Kashima}
\affil[1]{
Graduate School of Informatics, 
Kyoto University, Kyoto 606-8501, Japan}

\affil[*]{kashima@amp.i.kyoto-u.ac.jp}

\usepackage{bm}

\newcommand{\xini}{{{x}_0}} 
\newcommand{\tini}{0}
\newcommand{\xf}{\tilde{x}} 
\newcommand{\tf}{{\tau}}

\newcommand{\tL}{\tilde L}

\newcommand{\prj}{{\varrho}}
\newcommand{\Id}{\bmI}
\newcommand{\bmPi}{\Pi}

\newcommand{\osci}{\bmv}
\newcommand{\simu}[1]{({#1})}
\newcommand{\pseudo}{\tmk}

\newcommand{\xx}{{x}} 
\newcommand{\uu}{{u}}
\newcommand{\ff}{{f}} 
\newcommand{\bfg}{{g}} 

\newcommand{\bfe}{{e}}
\newcommand{\redst}{{z}}

\newcommand{\tmk}{^{T}}

\newcommand{\bme}{{e}}
\newcommand{\bmn}{{\xi}}
\newcommand{\bmv}{{\bf v}}
\newcommand{\bmA}{{A}}
\newcommand{\bmB}{{B}}
\newcommand{\bmG}{{G}}
\newcommand{\bmI}{{I}}

\newcommand{\calT}{{\mathcal T}}    
    
\newcommand{\calL}{{\mathcal L}}    

\newcommand{\Ee}{{\rm e}}   
\newcommand{\tr}{\mathop{\rm Trace}\nolimits}

\newcommand{\bbR}{{\mathbb R}}    
    
\newcommand{\R}{{\mathbb R}}    
\newcommand{\domain}{}

\newcommand{\exz}[1]{\left\langle{#1}\right\rangle}
\newcommand{\exu}[1]{\left\langle{#1}\right\rangle}

\newcommand{\Lt}{{L^{\tf}}}

\newcommand{\calLt}{{\calL^{\tf}}}

\newcommand{\calLtxt}{\calL^{\tf}(\xf)}
\newcommand{\noiselevel}{{T}}
\newcommand{\levels}{{\noiselevel_{\rm L}}}
\newcommand{\levell}{{\noiselevel_{\rm H}}}

\newcommand{\vmat}[2]{\left[
    \begin{array}{c}
      #1\\#2
    \end{array}\right]}
\newcommand{\dmat}[4]{\left[
    \begin{array}{cc}
      #1&#2\\
      #3&#4
    \end{array}\right]}

\newcommand{\hmat}[2]{\left[
    \begin{array}{cc}
      #1&#2
    \end{array}\right]}

\newcommand{\rev}[1]{#1}
\newcommand{\revsr}[1]{#1}
\newcommand{\ahmet}[1]{{#1}}
\newcommand{\template}[1]{}
\newcommand{\Letter}{{paper}}

\begin{abstract}
Control of nonlinear large-scale dynamical networks, e.g.,   
collective behavior of agents interacting via a scale-free connection topology, is a central problem in many scientific and engineering fields. For the linear version of this problem, the so-called controllability Gramian has played an \rev{important} role to quantify how effectively the dynamical states are reachable by a suitable driving input. In this \Letter, we first extend the notion of the controllability Gramian to nonlinear dynamics in terms of the Gibbs distribution. Next, we show that, when the networks are open to environmental noise, the newly defined Gramian is equal to the covariance matrix associated with randomly excited, but uncontrolled, dynamical state trajectories. This fact theoretically justifies a simple Monte Carlo simulation that can extract effectively controllable subdynamics in nonlinear complex networks. In addition, the result provides a novel insight into the relationship between controllability and statistical mechanics.
\template{(must be under 200 words)}
\end{abstract}
\begin{document}

\flushbottom
\maketitle
%
%
\thispagestyle{empty}

\template{\noindent Please note: Abbreviations should be introduced at the first mention in the main text – no abbreviations lists. Suggested structure of main text (not enforced) is provided below.}

\section*{Introduction}

Control, i.e., external forcing aimed at achieving desirable dynamical trajectories, of 
nonlinear large-scale dynamical networks is of 
major interest in many research fields such as gene regulatory networks, 
infection spreads, human brains, financial markets, smart grids, to list just a few \cite{Liu2011,Yuan2013,Power2011,Delpini2013}. 
To \revsr{investigate} how difficult such networks are to control, controllability, originally defined in control theory\cite{Kalman1963,Zhou1996}, has attracted much attention, mainly in physics research \cite{Liu2011,Yan2012,Yuan2013,Sun2013,Menichetti2014,Nepusz2012,Posfai2013,
Yuan2014,Zhao2015,Cowan2012,Wang2012}. 
Among them, \revsr{Kalman's controllabilty matrix has played an important role to determine whether every dynamical state of a linear system is reachable \cite{Liu2011}. Beyond this controllability \emph{determination}, }
the so-called controllability Gramian, that is only defined for linear systems, \revsr{provides} much of the \revsr{\emph{quantitative}} information concerning this problem. For example, every dynamical state is reachable if and only if the Gramian is nonsingular. Moreover, the minimum control energy required to drive the current state to a target one is represented as a quadratic form associated to the inverse of the Gramian, which is utilized to analyze the effect of connection topology \cite{Yan2012}. In this context, the condition number of the Gramian is a meaningful index to characterize the nonlocality of linear complex networks \cite{Sun2013}.

The controllability of complex networks with \emph{nonlinear} dynamics is also being actively investigated \cite{Whalen2015,Cornelius2013}. 
\revsr{The Lie bracket gives a natural extension of Kalman's controllability matrix rank condition for the controllability determination \cite{Whalen2015}.}
However, an analogous \revsr{controllability Gramian} for nonlinear dynamics has not yet been developed, even in control theory \cite{Antoulas2005,Schilders2008,Besselink2014,Astolfi2010,Scherpen1993,
Zuazua2014}, although the controllability Gramian of a linearized system is useful in some applications. 
\rev{One of only a few existing approaches is the empirical Gramian \cite{Lall2002} that appears in simulation-based model reduction methods \cite{Willcox2002,Kunisch2002} mainly developed in computational physics and numerical analysis. \revsr{The empirical Gramian is constructed using simulation data, which is in stark contrast to the controllability Gramian. Furthermore, it has been widely applied to \emph{nonlinear} large-scale systems} 
\cite{Lall2002,Hahn2003}.} However, although this is equal to the controllability Gramian when the dynamics are linear, there are no theoretical underpinnings for such an application to nonlinear cases.

The goal of this {\Letter} is to introduce a novel matrix measure for the controllability \revsr{quantification} of nonlinear network dynamics, to reveal its specific feature under stochastic noise, and to provide a simulation-based method for dynamical network reduction, together with its theoretical justification. 
To this end, we first extend the notion of the controllability Gramian to nonlinear systems \rev{from a statistical mechanics viewpoint}, and show the validity through its application to \revsr{controllability quantification}.  
Then we show that, when the network is open to environmental noise, the \rev{newly proposed} Gramian is equal to the covariance matrix of the uncontrolled dynamics. This equality brings about new insights into the relationship between controllability, simulation data, and stochasticity. 
This work is largely inspired by the path integral approach  proposed by Kappen \cite{Kappen2005}. 
Although this concept is not directly used as a numerical procedure to solve the optimal control problem below, it is a key building block to prove the main result. 

\template{
The Introduction section, of referenced text\cite{Figueredo:2009dg} expands on the background of the work (some overlap with the Abstract is acceptable). The introduction should not include subheadings.
}

\section*{Results}

\subsection*{Controllability function and Gramian.}
Consider the nonlinear controlled dynamics:
\begin{equation}\label{eq:controlled}
	\frac{d}{dt} \xx(t) = \ff(\xx(t))+\bfg(\xx(t))\uu(t),\ \xx(\tini)=\xini
\end{equation}
where $t$ represents time, $\xx(t)=[x_1(t),\ldots,x_n(t)]\tmk$ and $\uu(t)$ are the state 
and input variables, and smooth functions $\ff$ and $\bfg$ describe the autonomous dynamics and the input effect, respectively. \ahmet{In (\ref{eq:controlled}), the initial state $\xini$ is fixed, which affects both controllability determination and quantification. This can be  arbitrarily chosen, although typically the initial state is fixed to a stable equilibrium in the conventional controllability quantification results of nonlinear systems \cite{Scherpen1993}. 
Moreover, all the results in this {\Letter} hold for any probabilistic initial state (i.e., $\xini$ is a random variable) and multi input cases as far as $\xini$ is independent of the input noise below. }
For a final time $\tf>0$, the minimum control effort $\int_{\tini}^{\tf} \frac{1}{2}u^2(t) dt$ to achieve $\xx(\tf)=\xf$ is referred to as a controllability function denoted by $\Lt(\xf)$. 
When the dynamics are linear, i.e., $\ff(\xx)=\bmA\xx$ and $\bfg(\xx)=\bmB$ with constant matrices $\bmA,\ \bmB$ of compatible dimensions, the matrix $\bmG_\tf=\int_{\tini}^{\tf}\Ee^{\bmA s}\bmB \bmB\tmk\Ee^{\bmA\tmk s}ds$ is called the controllability Gramian. It is well known that $\Lt(\xf)=\frac{1}{2}\xf\tmk \bmG_\tf^{-1}\xf$ provided $\bmG_\tf$ is nonsingular when $\xini=0$ \cite{Antoulas2005,Sun2013,Zhou1996}. \rev{However, this definition of $\bmG_\tf$ cannot straightforwardly be extended to nonlinear systems.} Here, the controllability function and Gramian were introduced independently, and then a simple quadratic relation was shown. By changing our way of thinking, let us define $\bmG(\Lt)$, which we call Gibbs Gramian, in terms of the Gibbs distribution associated with the given controllability function
\begin{equation}\label{eq:Gramian}
	\bmG(L)= \int\domain \phi_L(\xf)\xf \xf\tmk d\xf,\ 	\phi_L(\xf)=\frac{\Ee^{-L(\xf)} }{\int\domain \Ee^{-L(\xf)} d\xf}.
\end{equation}
When $\Lt(\xf)=\frac{1}{2}\xf\tmk \bmG_\tf^{-1}\xf$, the Gaussian integral formula shows $\bmG(\Lt)=\bmG_\tf$. 
Therefore, this definition is consistent with the conventional one for linear dynamics. 
It should be noted that the definition of the controllability function does not assume linearity. 
Thus, this definition can readily be employed also for nonlinear cases. 
Another important feature is that we do not need to care about the reachability of each state.  
Even if some $\xf$ are not reachable by any finite energy input \rev{(e.g., linear dynamics for which $\bmG_\tf$ is singular), $\Lt(\xf)=+\infty$} causes no problem in (\ref{eq:Gramian})
because it simply leads to $\phi_{\Lt}(\xf)=0$. This means we can handle network dynamics that evolve on a specific domain due to the dynamics' structure or physical constraints. 

\ahmet{For linear systems with the initial state at the origin,  eigenstructure of the controllability Gramian $\bmG_\tf$ is useful for identifying directions in the state space that require small control energy to be reached. The Gibbs Gramian enjoys a similar property. Specifically,  principle component analysis on $\bmG(\Lt)$ reveals all effectively reachable directions. For instance, by setting $\Lt(\xf)=\frac{1}{2}\xf\tmk \bmG_\tf^{-1}\xf$, we observe for the linear case that the principle eigenvector of $\bmG(\Lt)=\bmG_\tf$ minimizes $\frac{L^{\tau}(\xf)}{\| \xf\| ^ 2}$, that is, the control effort divided by the squared distance takes its minimum value when the final state $\xf$ lies on the principle eigenvector. Another interpretation is that, of every possible direction, with a fixed control energy the state can be driven the furthest from the origin by driving it to a destination state that lies along the principle eigenvector. 
Interpretation for the nonlinear case has similarities with the linear case. 
Note that, for any unit vector $\bme$, large $\phi_\Lt(\xf) |\bme\tmk \xf|^2$ implies that a small energy input  can be used (i.e., small $\Lt(\tilde \xx)$, and consequently large $\phi_\Lt(\xf)$) to place the state far from the origin along the direction of $\bme$ (i.e., large $|\bme\tmk \tilde \xx|^2$) at the final time.  
Then, its spatial integral over all final states $\xf$ satisfies the following theorem, which readily follows from the equality $\bme\tmk\bmG(\Lt)\bme=\int\domain \phi_\Lt(\xf) |\bme\tmk \xf|^2 d\xf$:

\emph{The integral $\displaystyle \int\domain \phi_\Lt(\xf) |\bme\tmk \xf|^2 d\xf$ is maximized when $\bme$ is the principal eigenvector of $\bmG(\Lt)$. }

\noindent  In this sense, the principal eigenvector of the Gibbs Gramian captures the direction along which the states can be reached furthest from the origin by a control effort that is small on average. In addition, it is trivial to change the reference point. For example, one can modify the definition as $\bmG(L)= \int\domain \phi_L(\xf)(\xf-\xini) (\xf-\xini)\tmk d\xf$ to evaluate the travelling distance instead of the distance from the origin. Similarly, the secondary, and further, eigenvectors enable us to characterize an effectively reachable subspace. 
}
See also the subsection entitled \emph{Dimensionality reduction of nonlinear network dynamics} below for another interpretation in terms of the minimal projection error. The conclusion is that the Gibbs Gramian introduced $\bmG(\Lt)$ in (\ref{eq:Gramian}) is a proper extension of the conventional controllability Gramian $\bmG_\tf$ for nonlinear dynamics.

\subsection*{Stochasticity connects Gibbs Gramian and simulation data.}
\begin{figure}[t]
	\centering
	\includegraphics[width=13cm]{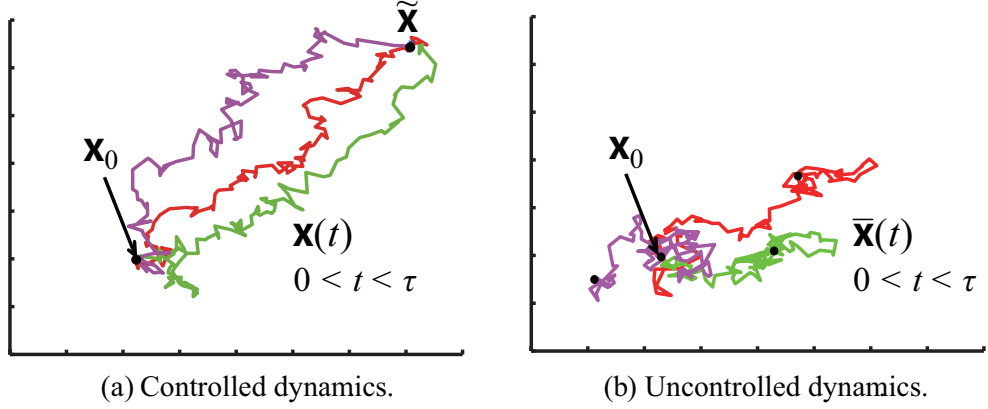}
	\caption{Typical behavior of controlled and uncontrolled dynamics open to environmental noise. (a) Sample paths of (\ref{eq:controlled_excited}) controlled by a fixed feedback law that regulates $\xx(\tf)=\xf$. The corresponding control effort is measured by $\exu{\int_0^\tf \frac{1}{2}u^2(t)dt}$, which is the average of $\int_0^\tf \frac{1}{2}u^2(t)dt$ over these sample paths. Then, $\calLtxt$ is the minimum of these average values over all such control laws. (b) Sample paths of the noise response $\bar\xx(t)$ in (\ref{eq:excited}).}
\label{fig:ctrl}
\end{figure}

For linear dynamics, $\bmG_\tf$ is given as a solution to a linear matrix equation.  
\ahmet{On the other hand, the controllability function $\Lt$ is given as a solution to a nonlinear partial differential equation for nonlinear dynamics\cite{Scherpen1993}. Therefore, }
it is not realistic to compute $\Lt$, and consequently $\bmG(\Lt)$, even for small-scale cases. 
However, the situation drastically changes when the input is disturbed by random noise:
\begin{equation}\label{eq:controlled_excited}
	\frac{d}{dt} \xx(t) = \ff(\xx(t))+\bfg(\xx(t)) (\uu(t)+ \sqrt{\noiselevel}\bmn(t)),\ \xx(\tini)=\xini
\end{equation}
where $\noiselevel>0$ is a noise level or temperature, $\bmn(t)$ is white noise such that $\exz{\xi^2(t)}=1$, and the expectation is taken over noise samples \cite{Karatzas1998,Masuda2010}.
\revsr{There are several theoretical results concerning the controllability determination of stochastic systems (e.g., approximate controllability \cite{Buckdahn2006,Mahmudov2003}). In this paper, we define a \emph{stochastic controllability function} $\calLtxt$ for the controllability quantification as 
\begin{align}\label{eq:stoc_ctrl_fcn2}
	& \calLtxt= \inf_{u} 
	\exu{ \int_{\tini}^\tf \frac{1}{2}\|\uu(t)\|^2 dt + \Phi(\xx(\tf)-\xf)},
\end{align}
where the infimum is taken over all feedback control laws and we define $\Phi$ such that 
\begin{equation}\label{eq:Phi}
	\Phi({0})=0,\ \Ee^{-\Phi(\xx)}\propto \delta(\xx),\ \xx\in \bbR^n
\end{equation}
with the Dirac's delta function $\delta$. Then,  
the terminal cost is an alternative representation of the terminal boundary constraint $\xx(\tf)=\xf$, since $\Phi(\xx(\tf)-\xf)=+\infty$ when $\xx(\tf)\neq\xf$. Therefore, $\calLtxt$ can be viewed as}
the minimum expected value of the control effort 
to \revsr{regulate} $\xx(\tf)=\xf$; see Figure \ref{fig:ctrl}a. \revsr{Similarly to the deterministic case, we do not require the boundedness of $\calLt$.} It should be emphasized that the resulting Gibbs distribution $\phi_\calLt$ is not identical to  $\phi_\Lt$, and depends on $T$. Next, we refer to the uncontrolled ($\uu(t)=0$), but randomly excited dynamics $\bar\xx(t)$ as noise response 
\begin{equation}\label{eq:excited}
	\frac{d}{dt} \bar \xx(t) = \ff(\bar \xx(t))+\bfg(\bar\xx(t))\sqrt{\noiselevel} \bmn(t),\ \bar\xx(\tini)=\xini,
\end{equation}
whose sample path is shown in Figure \ref{fig:ctrl}b. 
The key finding in this {\Letter} is \revsr{the following theorem, the proof of which is in the \emph{Methods} section:

\emph{The probability density function of $\bar\xx(\tf)$ is given by $\phi_{\calLt/T}(\bar\xx)$, that is, the noise response $\bar\xx(\tf)$ obeys the Gibbs distribution associated with $\calLt/T$.}}

This result means that the noise response data completely characterizes the minimum required input energy $\calLtxt$ for each target state $\xf$.
An intuitive reason for this nontrivial relation \rev{to hold} is that the noise in (\ref{eq:excited}) is added through the input channel. This type of noise is known to have an ability to search for the \rev{solution to a wide class of optimal control problems} \cite{Kappen2005}.
By this connection, the noise response data $\bar\xx(t)$ inherently contains information about the \rev{control energy minimization} problem. Therefore, this bridges the gap to the controllability function that is defined via the \rev{minimum energy control input}.  

Note that the evaluation of $\calLtxt$ over the whole state space based on the density function estimation of $\bar\xx(\tf)$ is still computationally intractable. However, in this {\Letter}, we focus on the Gramian \revsr{induced by the stochastic controllability function}, \rev{which is given by the spatial integral in (\ref{eq:Gramian}), and is much easier to determine than the pointwise evaluation of $\calLtxt$. Actually, \revsr{the theorem above yields the following equality for the \emph{stochastic Gibbs Gramian} $\bmG(\calLt/T)$}:}  
\begin{equation}\label{eq:main}
	\bmG(\calLt/T)=\left\langle \bar\xx(\tf) \bar \xx(\tf)\tmk \right\rangle.
\end{equation}
This equality tells us that the \revsr{stochastic} Gibbs Gramian can easily be calculated via Monte Carlo sampling of the uncontrolled dynamics open to environmental noise.
\rev{Furthermore, both this computation and also the principle component analysis of $\bmG(\calLt/T)$ are} efficiently implementable because various algorithms to achieve computational scalability exist for both Monte Carlo sampling (e.g., importance sampling) and matrix eigenvalue analysis. 
\revsr{Thus, the novel equality (\ref{eq:main}) characterized in this paper leads to the first numerically tractable procedure to find an effectively reachable subspace of large scale nonlinear dynamics, when they are open to environmental noise}, and is particularly useful for network dynamics for which only simulation algorithms, or time-series data collected in a noisy environment, are available. 

\subsection*{Dimensionality reduction of nonlinear network dynamics.}

\ahmet{The controllability quantification enables us to characterize subspaces that require a large control energy to be reached. By eliminating such subspaces, we can obtain a reduced order model, which is expected to well approximate the state trajectories as long as the input energy is not large. Actually, the dimensionality reduction of (mainly linear\cite{Antoulas2005,Ishizaki2014}) dynamical systems, which is helpful for understanding the hidden core mechanism, or to perform efficient numerical simulation, is an important application of the controllability quantification. In this section, we investigate two conceptually different nonlinear model reduction methods in the light of the Gibbs Gramian. }

\revsr{Let an integer $k(<n)$ be the desired order of the reduced model and define the set of $(n\times k)$-matrices $\bmPi_k=\{\prj: \prj\pseudo \prj = \Id \}$ where $\Id$ denotes the identity matrix.
Suppose some $\prj\in\bmPi_k$ satisfies 
\begin{equation}\label{eq:lossless}
	(\Id-\prj \prj\pseudo) \xx(t) \approx 0 {\rm\ for\ all\ } t.
\end{equation}
Then, the reduced state $\redst(t)=[z_1(t),\ldots,z_k(t)]\tmk=\prj\pseudo \xx(t)$
can approximately recover the original one by $\xx(t)\approx \prj \redst(t)$. Hence, we expect the Galerkin projection given as $d{\redst}(t)/dt  = \prj\pseudo \ff(\prj \redst(t))+\prj\pseudo \bfg (\prj \redst(t))\uu(t)$ to be a good reduced order model of dynamics (\ref{eq:controlled}). 
In what follows, we focus on the problem of finding such a $\prj$. 
}

In computational physics, the
Proper Orthogonal Decomposition (POD), or Karhunen-Loeve method, has a long history of intensive research \cite{Willcox2002,Kunisch2002}. This is a simulation-based model reduction
method, and is widely used for the simulation of nonlinear large-scale dynamical systems as found in computational
fluid dynamics and aerospace engineering. Suppose we replace the requirement (\ref{eq:lossless}) by $\sum_{\tf\rev{\in\calT}}\|(\Id-\prj \prj\pseudo)\xx(\tf)\|^2$
where the error is evaluated at multiple, given time instances $\tf\rev{\in\calT}$.
\revsr{(This optimization criterion is equivalent to the maximal singular value of $(\Id-\prj \prj\pseudo)[\xx(\tf_1),\xx(\tf_2),\ldots]$.)} This is the fundamental idea of the POD, and is referred to as the method of snapshots. 
If we need to approximate only the autonomous system $d\xx(t)/dt=\ff(\xx(t))$, 
this optimization is computationally tractable even for nonlinear large-scale dynamics, and the resulting Galerkin projection yields a satisfactory reduced model. However, when controlled dynamics (\ref{eq:controlled}) are of interest, 
we need to determine which input signal $\uu(t)$ should be injected when collecting snapshots, because we cannot simulate the trajectories corresponding to all possible input signals. 
\revsr{Many practically useful techniques, as well as theoretical analysis tools, for this have been developed; see \cite{Holmes2012} and references therein. On the other hand, from a controllability quantification viewpoint, it is also reasonable to find $\prj$ such that $(\Id-\prj \prj\pseudo) \tilde \xx \approx 0$ if $\tilde\xx$ is reachable with a small energy input $u(t)$. 
For this purpose, the Gramian-based model reduction for linear systems employs a $\prj$ that maximizes $\tr(\prj \bmG_\tau \prj\tmk)$. The Galerkin projection associated with this choice extracts effectively reachable subdynamics, in that the resulting projection eliminates a subspace on which $\Lt(\xx)=\frac{1}{2}\xx\tmk \bmG_\tau^{-1} \xx$ is large. However, \ahmet{as mentioned at the beginning of the previous section, for nonlinear dynamics it is unrealistic to compute the controllability function $\Lt$, which is no longer a quadratic form. This is the main reason why there have been no practical methods for the control-theoretic model reduction of general nonlinear large-scale systems \cite{Holmes2012,Besselink2014,Astolfi2010,Scherpen1993}. 
This limited applicability}
shows a clear contrast to the POD. 
There are many results that attempt to solve optimal control problems by the POD \cite{Kunisch2004,Kunisch2005,Hinze2005}.
However, the relation between the simulation-based and Gramian-based model reductions has not yet been fully understood. 

The remainder of this section is devoted to forming a theory-bridge to connect these two model reduction approaches that were developed independently for similar purposes.
}
%
%
%
%
\revsr{
Concerning the input selection for the POD, 
}
the impulse signals for the empirical Gramian \cite{Lall2002}, or the sinusoidal signals for the frequency domain POD \cite{Willcox2002,Astolfi2010}, may be suitable for linear systems.  
\revsr{Actually, the POD with these input signals is equivalent to the Gramian-based model reduction for linear systems \cite[Chapter 5]{Holmes2012}, \cite[Section 9.1]{Antoulas2005}. }However, although it is technically easy to inject the same inputs for nonlinear systems, 
there is no solid justification for their use. 
An interesting solution is to choose white noise $\xi(t)$ for the input signal, and minimize the snapshots' ensemble average of the squared projection error, that is, $\exz{\sum_{\tf\rev{\in\calT}} \|(\Id-\prj \prj\pseudo)\bar \xx(\tf)\|^2}$. Note that (\ref{eq:main}) leads to 
\begin{equation}\label{eq:POD_criterion}
	\exz{ \sum_{\tf\rev{\in\calT}} \|(\Id-\prj \prj\pseudo)\bar \xx(\tf)\|^2} = \sum_{\tf\rev{\in\calT}} \int\domain \phi_{\calLt/T}(\xf) \|(\Id-\prj \prj\pseudo)\xf \|^2 d\xf. 
\end{equation}
Therefore, \revsr{the approximation error of the noise response data is equal to the projection error weighted by the Gibbs distribution associated with the stochastic controllability function $\calLt$. Consequently, the POD evaluates the projection errors on the trajectories that are realizable by a small control effort, without computing any minimum energy input. 
\ahmet{In this sense, the POD with noise response data can be regarded as an easily implementable nonlinear model reduction method that explicitly takes the controllability into account.}

Note that (\ref{eq:POD_criterion}) is equal to $\tr(\sum_{\tf\in\calT}\bmG(\calLt/\noiselevel))	-\tr\left(\prj\pseudo\sum_{\tf\in\calT}\bmG(\calLt/\noiselevel)\prj\right)$. 
\ahmet{Therefore, the stochastic Gibbs Gramian-based reduction (the maximization of $\tr\left(\prj\pseudo\sum_{\tf\in\calT}\bmG(\calLt/\noiselevel)\prj\right)$) is equivalent to the best approximation of effectively reachable states (the minimization of the right-hand-side of (\ref{eq:POD_criterion})). This is another justification for the conclusion that the Gibbs Gramian is a proper extension of the conventional controllability Gramian. }
}

Furthermore, equality (\ref{eq:POD_criterion}) holds even for any nonlinear projection in the place of $\prj \prj\pseudo$, although its optimization is nontrivial.
\revsr{In the case of linear systems, observability, a dual concept of the controllability, is also investigated, and often referred to as the balanced POD \cite{Holmes2012}. Extensions in this direction are currently under investigation.}

\section*{Discussion}

The dynamics' nonlinearity makes the controllability sensitive to $T$.
\rev{This temperature dependency is discussed in this section. First,}
$\calLt/T$ in (\ref{eq:main}) and (\ref{eq:POD_criterion}) indicates that the input cost is inversely proportional to the noise level, 
that is, a less noisy (accurate) control channel is more expensive. In particular, 
as $\noiselevel\rightarrow 0$, the criterion $\exz{\sum_{\tf\rev{\in\calT}}  \|(\Id-\prj \prj\pseudo)\bar \xx(\tf)\|^2}$ evaluates the error at the snapshots located almost on the trajectory of the autonomous system $d\xx(t)/dt=\ff(\xx(t))$. 
Its interpretation from a controllability perspective is as follows: The input weight $\noiselevel^{-1}$ becomes unboundedly large, and consequently the states reachable with small control energy 
are limited to a small neighborhood around the autonomous trajectory (the noise is negligible because $\noiselevel\rightarrow 0$). 

\begin{figure}[t]
	\centering
		\includegraphics[width=17cm]{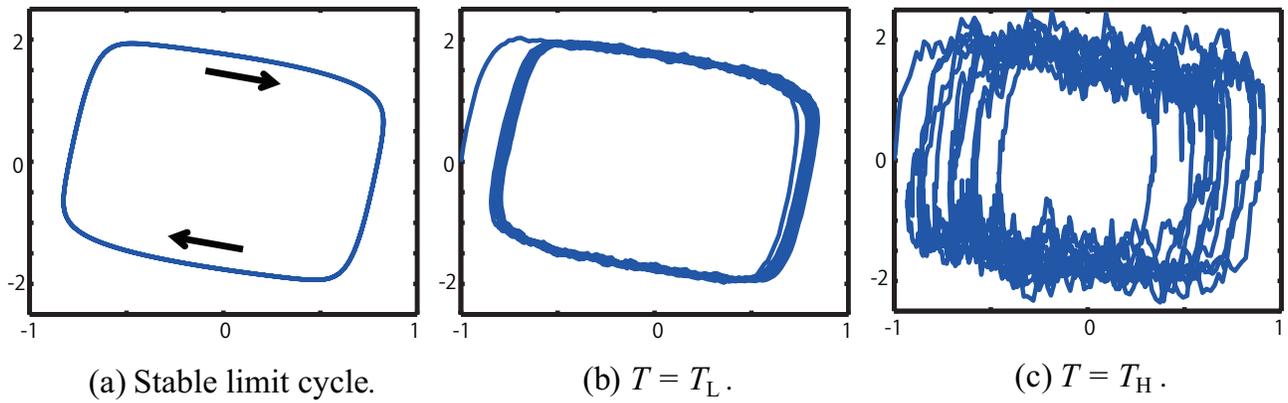}
\caption{A phase portrait of the FitzHugh-Nagumo neuronal oscillator in the $(v,w)$-plane. (a) The stable limit cycle of the noise-free individual dynamics. (b) A sample path of  $\osci_{1}(t)$ for $\noiselevel=\levels$. (c) A sample path of  $\osci_{1}(t)$ for $\noiselevel=\levell$. 
}
\label{fig:FH}
\end{figure}

\begin{figure}[t]
\centering
		\includegraphics[width=17cm]{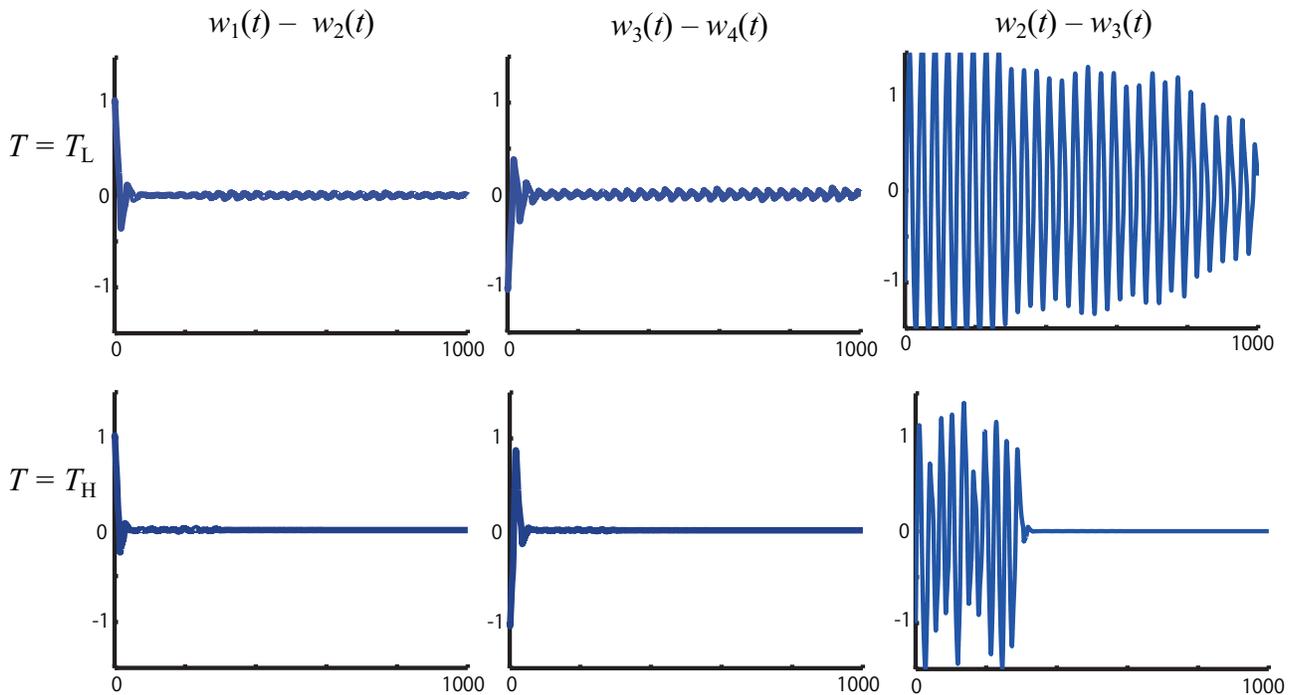}
\caption{(De)Synchronization phenomena in a sample path. 
 For both noise levels, $(w_1(t)-w_2(t))$ and $(w_3(t)-w_4(t))$ quickly decay due to the synchronization caused by the strong couplings. Synchronization is not observed in $(w_2(t)-w_3(t))$ for $T=\levels$ because the coupling strength $\eta_{23}$ is small. It shows a clear contrast to the quick noise induced synchronization for $T=\levell$. }
\label{fig:sync}
\end{figure}

Next, we demonstrate the nontrivial effect of the noise level by means of a numerical example of \revsr{$p$} identical, coupled neuronal oscillators of the FitzHugh-Nagumo model. The individual neuron generates the stable limit cycle shown in Figure \ref{fig:FH}a. 
The state variable of the $i$-th neuron is denoted by $\osci_i(t)=\hmat{v_i(t)}{w_i(t)}\tmk$, and the dimension of the entire system's state $\xx(t)=[\bmv_1(t)\tmk\ \bmv_2(t)\tmk\ \ldots\ \bmv_p(t)\tmk]\tmk$ is $n=2p$. The dynamics of the $i$-th neuron, subject to the diffusive coupling with nonuniform intensity and external forcing, are given by 
\begin{align*}
	& \frac{d}{dt}\vmat{v_i(t)}{w_i(t)} = \vmat{v_i-v_i^3/3-w_i}{0.08(v_i-0.8 w_i)}
	+\vmat{\sum\nolimits_{j\neq i} \eta_{ij}(v_j(t)-v_i(t))+\revsr{\hat u_i(t)}}{0} .
\end{align*}
\revsr{By using (\ref{eq:main}), we computed $\sum_{\tf\in\calT}\bmG(\calLt/\noiselevel)$ for $\rev{\calT=\{0.1,0.2,\ldots,1000\}}$ based on 1000 paths of the uncontrolled trajectories $\bar \xx(t)$ with $\uu(t)=0$, and its normalized eigenvector $\bfe_i$ corresponding to the $i$-th largest eigenvalue $\lambda_i$.  Low ($\levels=0.05^2$) and high ($\levell=0.5^2$) noise levels are considered. }
\revsr{Let $p=4$ and $\hat u_i=u(t)+\sqrt{\noiselevel} \xi(t)$ for all $i$, which means that only a common input is allowed.}
The symmetric coupling strengths $\eta_{ij}(=\eta_{ji})$ are given by $\eta_{12}=\eta_{34}=0.1,\ \eta_{23}=0.005$,
and $0$ for other pairs. The initial states are $\bmv_1(0)=-\bmv_3(0)=\hmat{-1}{0}\tmk,\ \bmv_2(0)=-\bmv_4(0)=\hmat{0}{2}\tmk$. 
For the uncontrolled trajectories $\bar \xx(t)$ with $\uu(t)=0$, 
apart from the fluctuation shown in Figure \ref{fig:FH}, we observed the following 3 (de)synchronization phenomena with a high probability: \simu{A}~$(\bmv_1-\bmv_2)$ and $(\bmv_3-\bmv_4)$ quickly decayed due to their strong couplings, \simu{B}~$(\bmv_2-\bmv_3)$ decayed only slowly for $\noiselevel=\levels$ because their coupling is weak, \simu{C}~$(\bmv_2-\bmv_3)$ quickly decayed for $\noiselevel=\levell$ because noise-induced synchronization occurred  \cite{Harada2010,Teramae2004}. See Figure \ref{fig:sync} for these phenomena observed in a sample path. 

As explained in Table \ref{tbl:correlation}, $\bfe_1$ and $\bfe_2$ approximately span the subspace given by
\begin{align*}
	& \noiselevel=\levels:\ \bmv_1 = \bmv_2=\dmat{0.0406}{-1.0199}{1.4383}{-0.4133} \bmv_3,\ \bmv_3 = \bmv_4, \\
	& \noiselevel=\levell:\ \bmv_1 = \bmv_2 = \bmv_3 = \bmv_4.
\end{align*} 
\revsr{Recall that $\prj=\hmat{\bfe_1}{\bfe_2}$ minimizes  
(\ref{eq:POD_criterion}) because $\tr\left(\prj\pseudo\sum_{\tf\in\calT}\bmG(\calLt/\noiselevel)\prj\right)$ is maximized; see the previous section}. Thus, the Galerkin projection onto this subspace extracts core subdynamics in the following two senses.
First, from a POD perspective, this subspace best approximates the noise response data; see the left-hand-side of (\ref{eq:POD_criterion}). This is confirmed by the fact that quick convergence to this subspace is nothing but the aforementioned (de)synchronization phenomena. Second, from a controllability perspective, this subspace best approximates the effectively reachable states; see the right-hand-side of (\ref{eq:POD_criterion}). In other words, even if we apply the optimal feedback control, it is expensive to avoid the (de)synchronization phenomena.
%

This can also be understood from the structure of the dynamics. 
Concerning \simu{A} and \simu{C}, since only the common input is allowed, the synchronization induced by the strong coupling and noise is difficult to prevent, independently of $\noiselevel$. On the other hand, concerning the desynchronization in \simu{B},
even though some well designed entrainment signals exist \cite{Harada2010}, 
they are not effective enough when the input weight $\noiselevel^{-1}$ is large. 
\revsr{As observed in
this example, controllability of highly nonlinear phenomena can be suitably captured from the noise-driven simulation data.
Note that reduced order models for controlled complex networks obtained by the proposed method do not always allow such a simple interpretation. In other words, this method can extract nontrivial core dynamics purely from time series data.}
%

In summary, we have proven that the noise response of the uncontrolled dynamics reveals the temperature-dependent controllability of general nonlinear network dynamics. 
\ahmet{This contribution consists of the following two achievements: One is a novel extension of the celebrated controllability Gramian for linear systems. To the author's best knowledge, this is the first nonlinear extension of the controllability Gramian, which was introduced over half a century ago and played a central role in the development of modern control theory \cite{Kalman1963}. The second achievement is equality (\ref{eq:main}), which mathematically proves that, when the system is open to environmental noise, the newly introduced Gramian is equal to the covariance matrix of the noise response data. This result forms a theory-bridge connecting controllability quantification and time series data analysis. An important outcome is that the equality (\ref{eq:POD_criterion}) yields an easily implementable method to a control-theoretic nonlinear model reduction problem for the first time. }
An extensive amount of noisy data is presently being gathered, and has been gathered to date, for a variety of uncontrolled systems. The equality (\ref{eq:main}) makes such data useful to glean insight into the modeling/controllability of controlled systems. 
We believe that this result can provide new methods and viewpoints in many research fields in view of the 
fact that much controllability related work is inspired by \rev{the pivotal contribution by Liu \emph{et al.} \cite{Liu2011}} For example, the condition number of the Gibbs Gramian should characterize the effect of nonlinearity on the network nonlocality analogously to the linear case 
\cite{Sun2013}. Also, 
in view of the numerical simulation above, the relation between the noise effect and the connection topology of the dynamical complex networks \cite{Ren2010} can be analyzed. 
Furthermore, 
\ahmet{the fact that any uncontrolled nonlinear dynamics subject to environmental noise obey the Gibbs distribution associated with $\calLt/T$, which is the minimum input energy divided by the temperature, suggests a nontrivial link to the \emph{canonical distribution} that is used in statistical mechanics}. 

\template{
The Discussion should be succinct and must not contain subheadings.
}

\section*{Methods}

\revsr{For simplicity of exposition we let $T=1$, but note that general results can be shown similarly.} 
It is well known \cite{Kappen2005} that the optimal value of the stochastic control problem in (\ref{eq:stoc_ctrl_fcn2}) satisfies
$\calLtxt=\tL(\tini,\xini)$ where the real scalar function $\tL(t,\xx)$ is the solution to 
the Hamilton-Jacobi-Bellman equation
\begin{align*}
	&\frac{\partial \tL}{\partial t}-\frac{1}{2} { \frac{\partial \tL}{\partial \xx}}\tmk \bfg \bfg\tmk \frac{\partial \tL}{\partial \xx}
	+\ff\tmk\frac{\partial \tL}{\partial \xx} 
	+\frac{1}{2}\tr 
	\left( \frac{\partial^2 \tL}{\partial \xx^2} \bfg \bfg\tmk \right)={0}, \\
	&\tL(\tau,\xx)=\Phi(\xx-\xf).
\end{align*}
Next, the logarithmic transformation $\psi(t,\xx) = \Ee^{- \tL(t,\xx)}$ yields 
the linear PDE
\begin{align*}
	& \frac{\partial \psi}{\partial t}+\ff\tmk \frac{\partial \psi}{\partial \xx}+\frac{1}{2}\tr 
	\left( \frac{\partial^2 \psi}{\partial \xx^2} \bfg \bfg\tmk \right) ={0}, \\
	& \psi(\tf,\xx)= \Ee^{-\Phi(\xx-\xf)}.
\end{align*}
This form allows us to apply the Feynman-Kac formula \cite{Karatzas1998} to obtain  
\begin{equation*}
	\psi(\tini,\xini) = \exz{ \Ee^{-\Phi(\bar\xx(\tf)-\xf)}}.
\end{equation*}
\revsr{Based on $\Ee^{-\calLtxt}=\Ee^{-\tL(\tini,\xini)}=\psi(\tini,\xini)$ and (\ref{eq:Phi}), 
we have 
\begin{equation}\label{eq:proof_phi}
	\phi_{\calLt}(\xf)
	=\frac{\exz{ \delta(\bar\xx(\tf)-\xf) }}{\int\domain \exz{ \delta(\bar\xx(\tf)-\xf) }d\xf}
=\frac{\exz{ \delta(\bar\xx(\tf)-\xf) }}{ \exz{ \int\domain\delta(\bar\xx(\tf)-\xf)d\xf }}
	=\exz{\delta(\bar\xx(\tf)-\xf)},
\end{equation} 
and consequently, 
\begin{equation}\label{eq:proof}
	\int\domain \phi_{\calLt}(\xf) w(\xf)d\xf = \int\domain \exz{\delta(\bar\xx(\tf)-\xf)} w(\xf) d\xf
	= \exz{ \int\domain \delta(\bar\xx(\tf)-\xf) w(\xf) d\xf}
= \exz{ w(\bar\xx(\tf))}
\end{equation} 
for an arbitrary smooth function $w(\xx)$ defined on $\R^n$, where we exchanged the order of expectation and spatial integral. 
The arbitrariness of $w(\xx)$ in (\ref{eq:proof}) means the probability density function of $\bar\xx(\tf)$ is given by $\phi_{\calLt}(\bar\xx)$. 
 Finally, (\ref{eq:proof}) with $w(\xx)=\xx\xx\tmk$ yields (\ref{eq:main}).}

\section*{Acknowledgements}

This work was in part supported by JSPS KAKENHI Grant Number 26289130, and by Japan Science and
Technology Agency, CREST.
\template{
Acknowledgements should be brief, and should not include thanks to anonymous referees and editors, or effusive comments. Grant or contribution numbers may be acknowledged.
}

\template{
\section*{Author contributions statement}

B.C. Hiesmayr has computed and written the results contained in this paper.

Must include all authors, identified by initials, for example:
A.A. conceived the experiment(s),  A.A. and B.A. conducted the experiment(s), C.A. and D.A. analysed the results.  All authors reviewed the manuscript. 
}

\section*{Additional information}

The author declares he has no competing financial interests.
\template{To include, in this order: \textbf{Accession codes} (where applicable); \textbf{Competing financial interests} (mandatory statement). 

The corresponding author is responsible for submitting a \href{http://www.nature.com/srep/policies/index.html#competing}{competing financial interests statement} on behalf of all authors of the paper. This statement must be included in the submitted article file.
}

%

%

\begin{table}[ht]
	\centering
  \begin{tabular}{|c | c | c |}
  \hline
     & $\noiselevel=\levels$ & $\noiselevel = \levell$ \\ \hline
    $\prj_1$ & $\dmat{-0.0805}{0.0629}{-0.4985}{0.4906}$ & $\dmat{-0.0619}{0.4947}{-0.4954}{-0.0697}$ \\ \hline
    $\prj_2$ & $\dmat{-0.0890}{0.0669}{-0.4869}{0.5013}$ & $\dmat{-0.0622}{0.4954}{-0.4966}{-0.0702}$ \\ \hline
    $\prj_3$ & $\dmat{0.0724}{0.0902}{0.4913}{0.4969}$ & $\dmat{-0.0611}{0.4968}{-0.4964}{-0.0533}$\\ \hline
    $\prj_4$ & $\dmat{0.0801}{0.0871}{0.4971}{0.4868}$ & $\dmat{-0.0611}{0.4976}{-0.4964}{-0.0531}$ \\
    \hline
  \end{tabular}
    \caption{\label{tbl:correlation}
The eigenvectors of $\sum_{\tf\in\calT}\bmG(\calLt/T)$ for the example. 
For both noise levels, $\lambda_i/\lambda_1<0.15$ for $i=3,4,\ldots,8$. The eigenvectors $\bme_1,\  \bme_2$ are given by $[\bme_1\ \bme_2]=[\prj_1\pseudo\ \prj_2\pseudo\ \prj_3\pseudo\ \prj_4\pseudo]\pseudo$ with $\prj_i$ listed above. 
Based on the standard correlation analysis, 
we conclude that $\prj_1\approx\prj_2,\ \prj_3\approx\prj_4,\ \prj_2\not\approx \prj_3$ for $\noiselevel=\levels$, and $\prj_1\approx\prj_2\approx\prj_3\approx\prj_4$ for $\noiselevel=\levell$.
}
\end{table}

\end{document}